\documentclass[reqno]{amsart}
\usepackage{amssymb}
\title{Stable manifolds of holomorphic diffeomorphisms}
\date{\today}
\author{Mattias Jonsson and Dror Varolin}
\address{Department of Mathematics; University of Michigan; Ann Arbor, MI
48109-1109}

\newcommand{\noi}{\noindent}
\newcommand{\rmk}{{\it Remark. }}
\newcommand{\bs}{\medskip}
\newcommand{\ie}{i.e.\ }
\newcommand{\eg}{e.g.\ }

\newcommand{\cg}{\mathcal{G}}

\newcommand{\cK}{\mathcal{K}}

\newcommand{\cN}{\mathcal{N}}
\newcommand{\co}{\mathcal{O}}

\newcommand{\cw}{\mathcal{W}}

\newcommand{\vp}{\varphi}

\newcommand{\Z}{\mathbf{Z}}
\newcommand{\R}{\mathbf{R}}

\newcommand{\C}{\mathbf{C}}
\newcommand{\N}{\mathbf{N}}
\newcommand{\p}{\mathbf{P}}

\newcommand{\cn}{\mathbf{C}^n}

\newcommand{\ck}{\mathbf{C}^k}

\newcommand{\relcomp}{\subset \subset}
\newcommand{\diff}{\mathrm{Diff}_{\co}}

\newcommand{\emb}{\hookrightarrow}

\newcommand{\klm}{\cK (\underline \lambda ,\underline m)}
\newcommand{\ul}{\underline\lambda}
\newcommand{\um}{\underline m}

\newtheorem{mthm}{Theorem}
\newtheorem{thm}{Theorem}[section]

\newtheorem{lem}[thm]{Lemma}
\newtheorem{defn}[thm]{Definition}
\newtheorem{prop}[thm]{Proposition}

\theoremstyle{definition}
\newtheorem*{remark*}{Remark}
\addtocounter{section}{-1}             % Start with section 0
\numberwithin{equation}{section}       % Number formulas within sections

\begin{document}

\thispagestyle{empty}

\maketitle

\section{Introduction}\label{intro}

Let $M$ be a complex manifold and fix once and for all a complete
Riemannian metric on $M$ and a holomorphic diffeomorphism, or
automorphism, $f \in \diff (M)$.  Recall that the 
\emph{stable manifold} 
$W^s_p$ through a point $p \in M$ with bounded orbit is
defined by
$$W^s_p := \left \{ x \in M\ |\ \mathrm{dist}(f^N x,f^N p) \le C \rho^N
\ \mathrm{for\ } N\ge0 \right \},$$ 
where $\rho=\rho_p<1$ and $C=C_p>0$. It
turns out that often $W^s_p$ is an immersed complex manifold.
Assuming this to be the case, the following problem was posed by
E. Bedford~\cite{bedford}.
\begin{center} 
  \textbf{Problem:} \textit{Determine the complex structure of
    the stable manifolds of $f$.}
\end{center}

In many cases it can be shown that $W^s_p$ is a monotone union of balls,
and this in turn implies~\cite{brown} that it is diffeomorphic to real 
Euclidean space. Moreover, by the contracting nature of the dynamics, 
one sees that the Kobayashi pseudometric of $W^s_p$ vanishes identically. 
However, when $\dim(W^s_p)\ge 3$, it is not possible to deduce only from these
properties that $W^s_p$ is biholomorphic to Euclidean space.
For example, there exist monotone unions of balls which are not Stein
\cite{fornaess}.  (The question of Steinness of monotone unions of
balls in complex dimension 2 is open.)  When dim$(W^s_p)=1$, the 
Uniformization Theorem implies that 
$W^s_p$ is biholomorphic to $\C$~\cite{bls,wu}.

The main results of this paper are proved in the non-uniform setting, \ie
with respect to compactly supported invariant measures.  More precisely,
we say that a subset $A \subset M$ is \emph{invariant} if $fA=A$, and that it
has \emph{total measure} if $\mu (A) = 1$ for every compactly supported
invariant probability measure $\mu$.

Our main objective in this paper is to prove the following theorem.

\begin{mthm}\label{main}
  There exists an invariant Borel set $\cK (f) \subset M$ of total
  measure such that for every $p \in \cK (f)$, $W^s_p$ is a complex
  manifold biholomorphic to complex Euclidean space.
\end{mthm}

Let us clarify things a little further.  The set $\cK (f)$ is the set of
so-called \emph{Oseledec points} or \emph{regular points}.  Its existence is
a part of the well known and fundamental theorem of V. Oseledec 
\cite{oseledec}.  The equally fundamental work of Pesin~\cite{pesin} states
in part that the stable manifold passing through each Oseledec point $p \in
\cK (f)$ is an immersed (complex) manifold.  What we show is that for each 
$p\in\cK(f)$, $W^s_p$ is biholomorphic to $\ck$, 
where $k=\mathrm{dim}_{\C}W^s_p$.

We postpone to Section~\ref{lyap} a more detailed discussion of the results of 
Oseledec and Pesin which we will use in this paper.  For now, however, we
content ourselves with saying that every point in $\cK (f)$ has a bounded 
orbit, although $\cK(f)$ itself need not be bounded.  Moreover, we emphasize
that, to have a non-trivial result, it is necessary to have at least one 
invariant measure with compact support, but this is guaranteed to happen
once $f$ leaves invariant a bounded subset of $M$.

Our approach to proving Theorem~\ref{main} is to associate to the
dynamical system $f$ a certain ``unraveled'' dynamical system, and
then conjugate the latter to a much simpler (polynomial) 
dynamical system on
the so-called stable distribution. To state the result more
precisely we need to develop some notation and concepts, which we
now proceed to do.

Recall that the stable distribution $E^s$ is a family of vector
subspaces $E^s_p$ of $TM_p$ on which $df$ is asymptotically
contracting: for $p\in M$ with bounded orbit, $E^s_p$ is given by
$$E^s_p := \left \{ v\in TM_p\ |\ |df^Nv|\le C\rho^N
\ \mathrm{for\ } N\ge0 \right\},$$ 
where $\rho=\rho_p<1$ and $C=C_p>0$.
Notice that this exponential decay is not uniform, 
\ie $p\mapsto\rho_p$, $p\mapsto C_p$ and even $p\mapsto\dim E^s_p$ 
could be discontinuous.  

From here on, let $\cK (f)$ denote the set of Oseledec points 
(whose well known definition is recalled in Section~\ref{lyap}).  
If the stable manifold $W^s_p$ defined
above exists as an immersed submanifold of $M$, then $(TW^s_p)_p=E^s_p$.  In
view of Pesin's work mentioned above, this is the case whenever $p\in \cK (f)$.

Next we discuss what was meant by "unraveled" above.  To this end,
even though it might happen that for $q \neq p$ the two stable
manifolds $W^s_p$ and $W^s_q$ intersect and thus agree, we treat
$W^s_p$ and $W^s_q$ as distinct stable manifolds. More precisely,
we define the set
\begin{equation}\label{Me6}
  \cw^s := \bigsqcup_{p \in \cK (f)} W^s_p.
\end{equation}
Since $W^s_{fp} = fW^s_p$, $f$ induces a bijection on $\cw^s$,
still denoted $f$, which is holomorphic on the fibers.  (We will
topologize $\cw^s$ and $E^s$ shortly.)  We note again that the bundle
$\cw^s$ is defined only over the topological space $\cK (f)$ of 
Oseledec points.  The following is our main result.

\begin{mthm}\label{pesin-u}
  There exists a measurable isomorphism $\Psi : \cw^s \to E^s|\cK (f)$ and
  a bundle automorphism $P : E^s|\cK (f) \to E^s| \cK (f)$ such that, for 
  every $p \in \cK (f)$, 
  \begin{enumerate}
  \item $P_p:E^s_p \to E^s_{fp}$ is a polynomial automorphism;
  \item $P^N_p \to 0$ locally uniformly on $E^s_p$;
  \item $\Psi_p:=\Psi|W^s_p$ is a biholomorphism of
    $W^s_p$ onto $E^s_p$ and $(d\Psi_p)_p=\mathrm{id}$;
  \item $\Psi \circ f \circ \Psi^{-1} = P$ on $\cK (f)$.
  \end{enumerate}
\end{mthm}

\begin{remark*}
  Note first that (3) implies Theorem~\ref{main}. 
  Secondly, we do not claim that $p\mapsto\mathrm{deg}(P_p)$ 
  is constant or even bounded. However, it is constant along
  orbits.  Finally, the map $\Psi$ turns out to be slightly better
  than measurable.  It is slowly varying; a concept we shall discuss
  more thoroughly in Section~\ref{slow}.
\end{remark*}

We now return to the question of topologizing $\cw^s$ and $E^s$.
Let
$$\mathrm{dist} (x,y) = \left \{
  \begin{array}{l@{\quad}r}
    \mathrm{dist}_p (x,y) & x,y \in W^s_p \\
    \mathrm{dist}_p (x,p)+ \mathrm{d} (p,q) + \mathrm{dist}_q (q,y) &
    x \in W^s_p,\ y \in W^s_q
  \end{array} \right .$$
$$\mathrm{dist} (v,w) = \left \{
  \begin{array}{l@{\quad}r}
    |v-w| & v,w \in E^s_p \\
    |v|+ \mathrm{d} (p,q) + |w| & v \in E^s_p,\ w \in E^s_q
  \end{array} \right .$$
The function $\mathrm{dist}_p$ is a distance on $W^s_p$ associated to the
complete Riemann metric on $M$ (so that $\mathrm{dist}_p$ 
recovers the intrinsic
topology of $W^s_p$). The function d appearing on the right
hand side of these definitions is the same for both functions. In
the case of Theorem~\ref{pesin-u}, we shall take d$(p,q)=\delta
_{p,q}$, \ie $\cw^s$ is the disjoint union of all the stable
manifolds.  We refer to this as the \emph{discrete case}.  Later
on, it is also useful to take d to be the Riemannian distance
on $M$.  We shall refer to this as the \emph{bouquet case}.  To
these topologies we associate the Borel sets, and it is with
respect to this $\sigma$-algebra that $\Psi $ is measurable.

We turn now to the hyperbolic picture.  
Recall that $f\in\diff(M)$ is \emph{hyperbolic on a compact set} 
$K$ if $K$ is invariant and there exists a continuous splitting
$$TM|_K=E^s\oplus E^u$$
with the following properties.
\begin{enumerate}
\item[(a)] $E^s_p$ and $E^u_p$ have constant rank for all $p\in K$, say 
  $k$ and $n-k$ ($n=\dim M$);
\item[(b)] $df (E^s_p) = E^s_{fp}$ and $df (E^u_p) = E^u_{fp}$ 
  for all $p\in K$;
\item[(c)] there exist positive constants $C$ and $\rho$, with $\rho <1$,
  so that, for all $p\in K$ and all $N\ge0$
  \begin{equation*}
    \Vert df^N|E^s_p\Vert\le C\rho^N
    \quad\mathrm{and}\quad
    \Vert df^N|E^u_p\Vert\ge C^{-1}\rho^{-N}.
  \end{equation*}
\end{enumerate}
In this case we write 
$$\cw^s=\bigsqcup_{p\in K}W^s_p \quad \mathrm{and}\quad E^s=\bigcup_{p\in K}
E^s_p.$$
Then $E^s$ is a continuous vector bundle and the set of (local) 
stable manifolds form a lamination near $K$ (see, \eg~\cite{shub}).
(We are abusing language slightly here: we are interested in the restriction
of $E^s$ to the set $K$ of hyperbolicity, but $E^s$ is being written instead
of $E^s |K$.)

Theorem~\ref{main} implies that 
$W^s_p$ is biholomorphic to $\C^k$ for every $p\in\cK(f)\cap K$. 
Moreover there are always invariant measures on an invariant 
compact set, so the set $\cK(f)\cap K$ is nonempty:
it contains every periodic point in $K$ and many more points
(unless $K$ is a finite set).
Nevertheless, one would like to prove that \emph{every} stable
manifold through $K$ is biholomorphic to $\ck$. We conjecture that
this is indeed the case. The main problem in proving this is that
even though $df^N$ uniformly contracts vectors in $E^s_p$, $p\in K$,
the exact rate of contraction can be highly nonconstant. In fact,
controlling the (asymptotic) rates of contraction within
$E^s$ is central to our approach.

On the other hand, the conjecture above is easy to prove if an 
(unfortunately quite strong) hypothesis is placed on the map $f$. 
If $A$ is the restriction to $E^s$ of $df$, let
\begin{equation*}
  L_+:=\limsup_{N\to\infty}\sup_{p \in K}N^{-1}\log\Vert A^N_p\Vert
  \qquad\text{and}\qquad
  L_-:=\liminf_{N\to \infty}\inf_{p\in K}-N^{-1}\log\Vert A^{-N}_p\Vert.
\end{equation*}
Note that $L_-\le L_+ <0$. 
We say that $f$ is \emph{equi-contracting} if $2 L_+ < L_-$.

Note, in particular, that if $f$ has one dimensional stable
manifolds, then $f$ is automatically equi-contracting.

\begin{mthm}\label{u}
  If $f$ is hyperbolic and equi-contracting on a compact set $K$, then there 
  exists a homeomorphism $\Psi : \cw^s  \to E^s$ over $K$ such that
  \begin{enumerate}
  \item for every $p \in K$, $\Psi_p:=\Psi|W^s_p$ is a biholomorphism of
    $W^s_p$ onto $E^s_p$ and $(d\Psi_p)_p=\mathrm{id}$;
  \item $\Psi \circ f \circ \Psi^{-1} = df|E^s.$
  \end{enumerate}
\end{mthm}
\noi In Theorem~\ref{u} we use the bouquet topology on $\cw^s$ and $E^s$. 
Note that the equi-contracting hypothesis implies that $f$ can be 
brought to a \emph{linear} form as opposed to the more general 
polynomial form given by Theorem~\ref{pesin-u}.

\bs

Regarding history, while Oseledec/Pesin theory has been used in complex
dynamics before (see \eg~\cite{bd},\cite{bls}), to our knowledge this
is the first application to the study of the complex structure of stable
manifolds in higher dimensions. 
When $f$ has a fixed point, the fact that $f$ is
conjugate to a normal form is due to S. Sternberg~\cite{st}, and
(independently, though much later) to Rosay and Rudin~\cite{rr} in
the holomorphic case with essentially the same proof. (We note that, 
even though Sternberg's theorem is stated for an attracting fixed point,
say $p$, one can reduce to this case  by restricting $f$ to its stable 
manifold $W^s_p$, since the latter is invariant in the fixed point case.)
The condition for linearization was known to C.~Siegel~\cite{siegel}.
In the general, non-stationary case, very little seems to have
been done. The main work we know of is due to M. Guysinsky and 
A.~Katok~\cite{g,gk}. However, they place rather strong hypotheses on
the spectrum of the of $df$ which, while sufficient (and perhaps
more so necessary) for their applications, would be much too
strong for the problem we are interested in here.

Roughly speaking, our approach combines the ideas from the proof of 
Sternberg's Theorem with techniques from Oseledec-Pesin Theory. 
The proof of Sternberg's theorem, as in~\cite{st} or~\cite{rr},
uses linear algebra to split the stable space $E^s_p$ into invariant
subspaces where $df$ has an essentially fixed rate of 
contraction (given by the eigenvalues of $df|E^s_p$). 
In the setting of Theorem~\ref{pesin-u} we use Oseledec-Pesin theory
in order to control the rate of contraction of $df$.

The organization of the paper is as follows.  In Section~\ref{exp}
we define a continuous family of uniformly sized charts for the
stable manifolds, and in Section~\ref{u-proof} we prove Theorem 3.
In Section~\ref{lyap} we state the results we use from the
Oseledec/Pesin theory and in Section~\ref{slow} we set up the {\it
slowly varying} formalism, working out some useful lemmas and
propositions.
%Section~\ref{exp2} contains the non-uniform analog
%of the result of Section~\ref{exp}.
In Section~\ref{loc-conj} we
prove the existence of $\Psi$ locally. This section is the main
step in the proof of Theorem~\ref{pesin-u}, the latter being
completed in Section~\ref{2pf}. \

 \bs

\noi {\bf Acknowledgments}  We thank John-Erik Forn\ae ss  and
Ralf Spatzier for their interest in this project and for interesting 
discussions, and more particularly we thank Ralf for also directing 
us to many useful references on related results in real dynamics.
We are also grateful to Charles Favre, Nessim Sibony and the 
referee for many helpful comments and suggestions.

\section{Holomorphic exponential maps I.  Hyperbolic case}\label{exp}
In this section, we construct a continuous (in $p$) family of
biholomorphic maps $\chi_p$ from a neighborhood of $0_p$ in
$E^s_p$ into $W^s_p$. To this end, let $f \in \diff (M)$ be
hyperbolic on a compact set $K \relcomp M$. For $\epsilon
> 0$, set
$$E^s_p (\epsilon) := \{ v \in E^s_p\ \left | \ |v|< \epsilon \} \right . \quad
\mathrm{and} \quad E^s (\epsilon) := \bigcup_{p \in K} E^s_p
(\epsilon), $$ the latter equipped with the bouquet
topology discussed in the Introduction.

\begin{prop}\label{isomorphism}
  There exists $\epsilon > 0$ and a continuous mapping $\chi : E^s ({\epsilon})
  \to \cw^s$ which maps each $E^s_p(\epsilon)$ biholomorphically into
  $W^s_p$, maps the zero vector $0_p \in E^s_p$ to $p\in W^s_p$, 
  and satisfies
  $d(\chi |E^s_p)_{0_p} = \mathrm{id}_{E^s_p}$.
\end{prop}

\noi \rmk  Originally, we had a rather complicated and not even completely
general proof of this proposition.  We thank C.~Favre for showing us a much 
simpler and complete proof, which we now present.
\begin{proof}[Proof of Proposition~\ref{isomorphism}.]  
  As mentioned before, it 
  is shown in~\cite{shub} that $\cw^s$ gives a lamination near $K$.  
  Cover $K$ by a finite number of balls $B_i,\ 1\le i\le l$. 
  Let 
  $$\xi_i:E^s(\epsilon_i)|B_i\cap K\to\cw^s,\quad1\le i\le l$$
  be local parameterizations of $\cw^s$ near $K \cap B_i$ such that, for all 
  $p \in B_i \cap K$, $\xi_i (0_p) = p$ and $d\xi_i (0_p) = id$.  Such 
  parameterizations exist if the balls $B_i$ are taken small enough.  
  Write $\epsilon := \min \{ \epsilon_i; 1 \le i \le l\}$.  Let $\cw^s 
 _i(\epsilon) = \xi_i(E^s
  (\epsilon) | B_i \cap K)$ be the image of 
  $E^s(\epsilon) | B_i \cap K$ 
  under $\xi_i$ , and set 
  $\cw^s(\epsilon )=\cup_i\cw^s_i(\epsilon)$.  
  Finally, let $\{ \vp_i ;\ 1 \le i \le l\}$ 
  be a partition of unity subordinate
  to the covering $\{ B_i ;\ 1 \le i \le l\}$ of $K$.  
  Then the map $\chi : E^s (\epsilon) \to \cw^s
  (\epsilon)$ whose inverse is given by the formula 
  $$\chi^{-1}(x)=\sum_{1\le i\le l}\vp_i(p)\xi_i^{-1}(x),\qquad x\in W^s_p,$$
  has the desired properties.
\end{proof}

\section{Proof of Theorem~\ref{u}}\label{u-proof}
We first sketch the basic idea of the proof.  Let 
$$A:=df|E^s$$ 
denote the restriction to $E^s$ of $df$ and, with $\chi$
as in Proposition~\ref{isomorphism}, set 
$$F:= \chi^{-1}\circ f\circ \chi \quad \mathrm{and} \quad \cw^s(\epsilon):= 
\chi (E^s (\epsilon)).$$  
Here and below, in order to avoid referring to the specific coordinates $v$ 
chosen on $E^s$, we use the notation $O(m)$ in place of the more common 
$O(|v|^m)$.

We want to show that the map 
$A^{-N}\chi^{-1} f^N$ converges, locally uniformly on 
$\cw^s$ as $N \to \infty$, to a biholomorphic map.  (Note that, because of
the use of $\chi^{-1}$, the former map is only defined on some 
compact subset of
$\cw^s$.)  Since $A^{-1} F - \mathrm{id} = O(2)$, we have, on a given compact
subset of $\cw^s$ and for $N$ sufficiently large, that
\begin{align*}
  A^{-(N+1)}\chi^{-1}f^{N+1}-A^{-N}\chi^{-1}f^N 
  &=A^{-N}\left(A^{-1}F-\mathrm{id}\right)\chi^{-1}f^N\\
  &\sim e^{(-L_-+2L_+)N},
\end{align*}
where the last estimate is uniform on compact sets.  
By the equi-contracting hypothesis, this implies locally uniform convergence. 
Injectivity and surjectivity of the limit map are then easily 
established.  The details are as follows.

Given $\delta>0$ with $L_++\delta<0$ 
there exists $N_0\in\N$ and $\epsilon>0$ such that
%then there exists $m\in \N$ such that for all
%$N\ge m$, $|A^N|\le e^{(L_++\delta)N}.$  One can choose $m$ so
%that, with $\epsilon > 0$ small enough, one has
$$|F^{N_0} v| \le e^{(L_+ +\delta)N_0}|v|\ \mathrm{whenever}\ v \in E^s (\epsilon).$$
This follows from the definition of $L_+$ and the fact that 
$(dF_p)_p=A_p$.
Now set
$$C:= \sup \left \{ \ |F^jv|/|v|\ ;\ 0\le j < N_0, v \in E^s (\epsilon)\ \right
\}.$$
For $N\ge0$, write $N=kN_0+j$ with $0\le j<N_0$. Then
$$|F^Nv| = |F^j(F^{kN_0}v)|\le C|F^{kN_0}v| \le Ce^{(L_++\delta)kN_0}|v|,$$
and so there exists $N_1=N_1(\delta)$ such that for all $N\ge N_1$,
$$|F^Nv|\le e^{(L_++2\delta) N},\quad v\in E^s(\epsilon) $$
Now consider a compact $J\subset\subset \cw^s$. By the contracting nature of
$f$ there exists $n\ge0$ such that $f^n(J)\subset\cw^s(\epsilon)$.  
Since $\chi^{-1}$ is continuous on $\cw^s(\epsilon)$, 
the above estimate implies that there exists 
$N_2=N_2(J,\delta)\ge N_1+n$ such that for all $N\ge N_2$,
$$|\chi^{-1}f^Nz|\le e^{(L_++3\delta) N},\quad z\in J $$

Since $(dF_p)_p=A_p$ there exists $C>0$ 
such that for all $v \in E^s(\epsilon)$,
$$\left | v-A^{-1} F v \right | \le C|v|^2.$$ 
Using the definition of $L_-$ and increasing $N_2$ if necessary, 
one then obtains, for all $z \in J$ and $N \ge N_2$, the estimate
\begin{align*}
  \left|A^{-N}\chi^{-1}f^Nz-A^{-(N+1)}\chi^{-1}f^{N+1}z\right|
  &\le\Vert A^{-N}\Vert\cdot\left|w_N-A^{-1}Fw_N\right|^2\\
  &\le e^{-(L_--\delta)N}C\left|w_N\right|^2\\
  &\le Ce^{(2L_+-L_-+5\delta)N},
\end{align*}
where $w_N = \chi^{-1} f^N z$. Since $2L_+<L_-$ by the equi-contracting
hypothesis, it follows that 
$$\Psi := \lim_{N\to \infty}  A^{-N} \chi^{-1} f^N$$
exists locally uniformly in the bouquet topology. Thus $\Psi$ is
continuous and clearly satisfies $d\Psi=\mathrm{id}$ as well as
the functional equation
$$\Psi\circ f\circ\Psi^{-1}=A.$$
We claim that $\Psi$ is in fact a homeomorphism. Clearly
$\Psi (W^s_p) \subset E^s_p$ for all $p\in K$. 
For fixed $p$, $\Psi |W^s_p$, being a uniform limit of automorphisms, 
is holomorphic and injective.  Thus $\Psi$ is itself injective. This implies
that $\Psi(\cw^s)$ contains a neighborhood $\cN$ of the zero section of
$E^s$. We now use the contracting property of $A$ to show surjectivity
of $\Psi$. Consider any $v \in E^s_p$ and pick $N$ large enough 
so that $A^Nv\subset\cN$, \ie there exists $y \in W^s_{f^Np}$ with 
$\Psi (y)= A^Nv$. Let $x := f^{-N} y$.  Then 
$$v = A^{-N} \Psi (y) = \Psi (f^{-N}y) = \Psi (x).$$ 
Thus $\Psi$ is a homeomorphism, which completes the proof.\qed

\section{Lyapunov data and stable manifolds}\label{lyap}
In this section, we give an overview, containing no proofs, of various results
in smooth ergodic theory.  There are several references which the reader can
consult for details.  We have taken most of our statements from 
\cite{pugh-shub}, but a more detailed proof of some of the theorems can be 
found in~\cite{mane}.

First, to an automorphism $f$ one can associate its Lyapunov data:  these are
vector spaces $E^{\lambda }_p \subset TM_p$, called the \emph{Lyapunov spaces}
of $f$, defined by
$$E^{\lambda }_p := \left \{ v \in TM_p \ \left | \lim_{N \to \pm \infty}
N^{-1} \log \left | df^N v \right | = \lambda \right . \right \}.$$
The numbers $\lambda = \lambda (p)$ such that $E^{\lambda}_p \neq \{ 0\}$ are
called the \emph{Lyapunov exponents}.

In general, of course, there are only a finite number of Lyapunov
exponents at a given point $p$.  A point $p$ with bounded orbit
such that
\begin{equation}\label{lsplit}
  TM_p = \bigoplus_{\lambda \in \R} E^{\lambda}_p,
\end{equation}
is called an \emph{Oseledec point} (or regular point) of $f$.  The
splitting (\ref{lsplit})  is called the \emph{Lyapunov splitting}.
We denote the set of Oseledec points by $\cK (f)$.

Before stating the basic result on Oseledec points, we need the
following definitions.

\begin{defn}\label{svf}
  Let $X \subset M$ be an $f$-invariant  Borel set.
  \begin{enumerate}
  \item
    A Borel function $R:X\to(0,\infty)$ is called 
    $\epsilon$-slowly varying if 
    $e^{-\epsilon}\le R(fp)/R(p)\le e^{\epsilon}$ for every $p\in X$.
 \item
    A collection of Borel functions 
    $\{R_{\epsilon}:X\to(0,\infty)\ ;\ \epsilon>0\}$ 
    is called a slow variation if $R_\epsilon$ is 
    $\epsilon$-slowly varying for every $\epsilon>0$.
  \item  A function $h:X\to(0,\infty)$ is called 
    slowly varying if there exists a slow
    variation $R_{\epsilon}$ such that 
    either $h \le R_{\epsilon}$ or $h \ge 1/R_{\epsilon}$ 
    for all $\epsilon>0$.
  \end{enumerate}
\end{defn}

\noi \rmk  In what follows, we shall have to control either the
growth or shrinking of certain functions along orbits of $f$.  To
distinguish these two situations, we establish the following
convention:  in the former case, the functions shall have ranges
of the form $(a,\infty)$ with $a \ge 1$, and in the latter, ranges
of the form $(0,b)$ with $b < \infty$.
\begin{thm}[\cite{oseledec}]\label{T2}
  The set $\cK (f)$ is an
  invariant Borel set of total measure.  Moreover, there is a slow
  variation $\{R_{\epsilon} : \cK (f) \to (1,\infty )\}$ such that for
  all $p\in\cK (f)$ and all $\epsilon>0$,
  \begin{itemize}
  \item[a)]
    \begin{equation}\label{Me9}
      R_{\epsilon}(p)^{-1}e^{-\epsilon N} 
      \le\frac{|d(f^N)_p v|}{e^{\lambda N}|v|}
      \le R_{\epsilon}(p)e^{\epsilon N}
      \quad\text{whenever}\quad
      v\in E^{\lambda}_p
    \end{equation}
  \item[b)]
    $$\measuredangle \left ( E^{\lambda}_p, E^{\lambda '}_p \right ) 
    \ge \frac{1}{R_{\epsilon}(p)}\ {whenever}\ \lambda ' \neq \lambda.$$
  \end{itemize}
\end{thm}

At every Oseledec point, one has the following decomposition.
$$TM_p = E^s_p\oplus E^0_p \oplus E^u_p,$$
where $$E^s_p = \bigoplus_{\lambda < 0} E^{\lambda}_p \qquad \mathrm{and}
\qquad E^u_p = \bigoplus_{\lambda >0} E^{\lambda}_p.$$
Given such a $p$ we define the \emph{stable manifold at $p$} by
\begin{equation*}
  W^s_p:= 
  \left\{x\in M\ \bigg| \ 
    \limsup_{N\to\infty}\frac1N\log\mathrm{dist}(f^Nx,f^Np)<0
  \right\}.
\end{equation*}
The Pesin stable manifold theorem can thus be stated as follows.
\begin{thm}[\cite{pesin}]
  For every $p\in \cK (f)$,
  $W^s_p$ is an immersed (complex) submanifold of $M$.  
\end{thm}
In fact, Pesin's result also tells us how the stable manifolds $W^s_p$
depend on the base point $p$. In particular we have a non-uniform
version of the exponential map in Proposition~\ref{isomorphism}.
\begin{thm}[\cite{pesin}]\label{T1}
  The stable lamination $\cw^s$ defined by~\eqref{Me6} is a
  slowly varying lamination on $\cK (f)$ in the following sense: let
  $M \emb \R^n$ be an isometric immersion into Euclidean space.
  Then there is a slowly varying function $r: \cK (f) \to (0,1)$
  with the property that if $D_p$ is the ball of radius $r(p)$  and
  center $0_p$ in $E^s_p$, then the orthogonal projection $\Pi :
  W^s_p \to E^s_p$  is invertible on the branch of $\Pi^{-1}
  (D_p)$ containing $p$. Moreover, there exists a map
  $\chi : E^s (r) \to \cw^s$ which maps $E^s_p(r(p))$ 
  biholomorphically into $W^s_p$, 
  maps $0_p$ to $p$, and satisfies $d\chi_{0_p}=\mathrm{id}_{E^s_p}$.
  Moreover, $p\mapsto\Vert d\chi_p|E^s(r(p))\Vert$ 
  is a slowly varying function on $\cK(f)$.
\end{thm}
\begin{remark*}
  This \emph{lamination} aspect of Pesin's theorem is
  rarely stated, but it is easily seen to be true if one 
  follows the proof, say in~\cite{pugh-shub}, based on 
  the graph transform.
\end{remark*}

The set $\cK (f)$ can be further decomposed into invariant subsets as
follows. For $l \in \N_+$, $\ul=(\lambda_1,\dots,\lambda_l)$
with $\lambda_l<\dots<\lambda_1<0$ and $\um = (m_1,\dots,m_l)$,
let
$$\klm := \left \{ p \in \cK (f)\ \left | \ E^s_p = E^{\lambda_1}_p \oplus
\dots\oplus E^{\lambda_l}_p \ \mathrm{and}\ \mathrm{dim}(E^{\lambda_j}
_p )= m_j\ 1 \le j\le l. \right . \right \}.$$ Then $$\cK (f) =
\bigcup_{\ul , \um} \klm,$$ and each $\klm$ is a Borel set which is 
invariant for $f$.  These subset of ``constant stable Lyapunov data''
will be crucial to our further analysis.

\section{Slowly varying bundles and maps}\label{slow}
In this section we establish definitions and basic results
about slowly varying objects.  The slowly varying
notion of regularity is the strongest form of 
regularity that can be expected to hold in the non-uniform picture. 
Roughly speaking, slowly varying objects can be treated as constants,
as long as we are interested in exponential estimates.

\subsection*{Measurable bundles and maps}
Recall that the relative $k$-Grassmannian of $TM$ is a bundle
$\cg_k (TM) \to M$ whose fiber over $p \in M$ is the set of $k$
dimensional complex subspaces of $TM_p$. A 
\emph{measurable complex vector bundle} 
over a Borel subset $X \subset M$ is then a
measurable section $E$ of the Grassmann bundle $\cg_k (TM)$ over
$X$. A measurable subbundle $E'$ of a measurable complex vector
bundle $E$ is a measurable complex vector bundle such that $E'_p$
is a subspace of $E_p$ for each $p\in X$.  We implicitly assume
that the base $X$ is invariant for $f$, and endow all such vector
bundles with the metric inherited from $TM$. The total spaces
are given the discrete topology discussed in the introduction.

As a matter of notation, given a function $g : X \to (0,\infty)$, let
$$E(g) := \bigcup_{p\in X} E_p (g(p)),$$
where $E_p(r) = \{ v \in E_p\ |\ |v|<r \}.$ 
Notice that this is a neighborhood of the zero section $O_E$
of $E$ and we call such a neighborhood a \emph{tube}.
Later on we will work with slowly varying tubes, \ie tubes
defined by slowly varying functions $g$.

Let $E\to X$ and $E'\to X$ be two measurable vector bundles.  
A (measurable) bundle map $\Phi : E\to E'$ 
\emph{fibered over $\phi : X\to X$} is then
a map of the total spaces, defined in some tube $E(g)$,
such that $\phi$ is measurable (with respect to the Borel
$\sigma$-algebra inhereted from $M$)
and $\Phi (E_p) \subset E'_{\phi(p)}$ for $p\in X$.
We shall require that $\Phi$ map the zero section of $E$ to the
zero section of $E'$, and that 
the map $\Phi_p := \Phi |E'_p$ be holomorphic for
each $p \in X$.  This situation is sometimes denoted by the shorthand 
$\Phi/\phi$.  
In our setting, it will always be the case that
$\phi=f$ or $\phi=\mathrm{id}$.

A bundle map $\Phi/\phi$ is said to be \emph{homogeneous of degree $m$} if
for every $c\in \C$ one has 
$\Phi (c \cdot v) = c^m\cdot \Phi(v).$ 
We use the notation
\begin{equation}\label{Me7}
  \Vert\Phi_p\Vert:= \max_{v \in E_p}\frac{|\Phi_p(v)|}{|v|^m},
\end{equation}
and the notation $\Vert\Phi\Vert: p \mapsto\Vert\Phi_p\Vert$.

More generally, $\Phi$ is said to be \emph{polynomial of degree $m$} if there
exist homogeneous maps $\Phi_j$ of degree $j$, $1 \le j \le m$
such that $\Phi = \sum_{j=1}^m \Phi_j$. By our requirement that
the fiber maps be holomorphic, every bundle map $\Phi/\phi$
has a homogeneous expansion 
\begin{equation}\label{Me8}
  \Phi=\sum_{m=1}^{\infty}\Phi_m
\end{equation}
Using the notation in~\eqref{Me7} we have 
$B(p):=\sup_m\Vert\Phi_{m,p}\Vert^{1/m}<\infty$ for every $p\in X$.
In general one can say very little about the dependence of $B(p)$ on $p$, 
but the maps we will work with have more regularity: we say that 
$\Phi$ is \emph{slowly varying} if $B$ is slowly varying in the sense of
Definition~\ref{svf}.

It trivially follows that if $\Phi$ is slowly varying, then so are
all of its homogeneous parts $\Phi_m$. Also, sums, 
compositions and inverses of slowly varying maps are easily seen
to be slowly varying.

\subsection*{Tubes associated to slowly varying maps}
Consider a slowly varying selfmap $T/\mathrm{id}$ 
of a bundle $E\to X$ such that
$dT|O_E=\mathrm{id}$. This has the homogeneous expansion
\begin{equation*}
  T=\mathrm{id}+\sum_{j=2}^{\infty}T_j.
\end{equation*}
Then $\Vert T_j\Vert\le B^j$ for some slowly varying function 
$B:X\to(2,\infty)$ and $T$ is defined on the slowly varying tube 
$E(1/B)$.

In the next two propositions, which will be crucial in the final
step of the proof of Theorem~\ref{pesin-u}, we will show that the range,
as well as the domain of injectivity of $T$ contain slowly varying tubes.
\begin{prop}\label{image-tube}
  Let $T/\mathrm{id}$ be as above.  Then there exists a slowly varying function
  $g:X\to(0,1)$ such that 
  $$T(E(1/B))) \supset E(g).$$
\end{prop}
\begin{proof} For each $v\in E(1/B)$, one has
  $$|Tv| \ge |v|-\sum_{j=2}^{\infty} (B|v|)^j = |v| - \frac{B^2|v|^2}{1-B|v|}
  =|v|\left ( 1 - \frac{B^2|v|}{1-B|v|} \right ).$$ Letting $|v| \le
  1/(2B^2) < 1/4B$, we see that $|Tv| \ge |v|/3$.  The proposition
  now follows (with $g=1/3B$) from this and the openness of the maps $T_p$. 
\end{proof}

Of course, since $dT_p=\mathrm{id}$, the Inverse Function 
Theorem says that $T_p$ is invertible on a neighborhood of $0_p$.
The next proposition shows that if $T$ is slowly varying then so is the size
of this neighborhood.
\begin{prop}\label{inj-tube}
  Let $T/\mathrm{id}$ be as above.  Then there exists a slowly varying
  function $h:X\to(0,1)$ such that $T_p$ is well-defined and
  injective on $E_p(h(p))$ for each $p\in X$.
\end{prop}
\begin{proof} 
  Let $h$ be a function to be specified shortly.  It is easily shown
  (say, using linear coordinates on $E_p$) that there is a constant $C$,
  depending only on the rank of $E$, such that for all $j\ge 2$, 
  $p\in X$, and $x,y \in E_p(h(p))$,
  $$| T_j (x) - T_j (y) | \le C (j+1)^kB(p)^j h(p)^{j-1} |x-y|,$$
  where $k = \mathrm{rank}(E)$.  Now
  $$\left | T(x) - T(y) \right | = \left | x-y+\sum_{j=2}^{\infty} (T_j(x) -
    T_j(y)) \right | \ge |x-y|\left ( 1 - \sum_{j=2}^{\infty} C (j+1)^k B
    (Bh)^{j-1} \right ).$$ Thus $h=1/(2B^3)$ does the trick provided
  $B$ is bounded from below by a sufficiently large constant. This
  completes the proof. 
\end{proof}

\subsection*{Contracting linear maps}
We say that a linear bundle map $A/\phi$ of a bundle $E\to X$ 
is \emph{contracting} if there exists $\lambda<0$ such that
$\limsup_{N\to\infty}\frac1N\log|A^Nv|\le\lambda$ for every
$v\in E$. If $A$ is slowly varying, then this implies that there
is a slow variation $R_\epsilon:X\to(1,\infty)$ such that
\begin{equation}\label{Me25}
  \Vert A^N_p\Vert\le R_\epsilon(p)e^{(\lambda+\epsilon)N}
  \quad p\in X,\ N\ge 1.
\end{equation}
There is a standard 
way of making this contraction more uniform by changing
the metric. To this end, fix $\epsilon>0$ and set
\begin{equation}\label{adapt1}
  \left\langle v,w\right\rangle^* 
  :=\sum_{N=0}^{\infty} e^{-2(\lambda-\epsilon)N}
  \left\langle A^Nv,A^Nw \right\rangle, 
  \quad v,w\in E
\end{equation}
and denote the associated norm by $|\cdot|^*$ and 
operator norm by $\Vert\cdot\Vert^*$.
Using~\eqref{Me25} it is straightforward to verify that
the series~\eqref{adapt1} converges, and that the metric thus
obtained is a Borel metric on $E$ with the following properties:
$\Vert A\Vert^* \le e^{\lambda+\epsilon}$ and 
there exists an $\epsilon$-slowly varying function 
$C:X\to(1,\infty)$ such that
\begin{equation*}
  |v|\le |v|^*\le C(p)|v|,\quad v\in E_p.
\end{equation*}

\subsection*{Splittings and flags}
Suppose that the vector bundle $E\to X$ \emph{splits}, \ie 
there exist subbundles $E_1, \dots , E_l$ of $E$
such that 
$$E = \bigoplus_{j=1}^l E_j.$$
We say that this splitting is \emph{slowly varying}
if the projection maps $E\to E_j$ are all slowly varying.
It is possible to show that 
the splitting $E=\oplus E_j$ is slowly varying if and only if the
angle functions $\measuredangle(E_i,E_j)$ are slowly varying functions.

To a splitting $\oplus E_i$ there is an associated \emph{flag} 
$V_{\bullet}$, \ie a sequence 
$V_{0}=\{0\}\subsetneq V_1\subsetneq\dots\subsetneq V_l=E$
of vector subbundles defined by 
$$V_j=\bigoplus_{i=1}^jE_i.$$
Since this flag comes from a splitting, it is equipped with
projections $\mathrm{pr}_j:E\to V_j$. We remark that the flag
$V_\bullet$ depends on the ordering of the vector spaces $E_i$.
For our applications, the $E_i$ are Lyapunov spaces and the Lyapunov
exponents provide a natural ordering.

\subsection*{Pseudo-linear maps}
Let $E\to X$ be a vector bundle with a splitting $E=\oplus E_j$.
A homogeneous selfmap $\Phi$ of degree $m$ of $E$ can then
be further decomposed as
\begin{equation}\label{flag-decomp}
  \Phi = \sum_{j=1}^{l} \sum_{|\alpha|=m}\Phi_{j,\alpha}
\end{equation}
where
\begin{equation}\label{Me23}
  \Phi_{j,\alpha}:E\to E_j
  \quad\text{and}\quad
  \Phi_{j,\alpha}(c\cdot v)=c^{\alpha}\Phi_{j,\alpha}(v).
\end{equation}
Here $c=(c_1,\dots,c_l)\in\C^l$, $\alpha=(\alpha_1,\dots,\alpha_l)\in\N^l$,
$c^\alpha=c_1^{\alpha_1}\dots c_l^{\alpha_l}$ and% acts on $\oplus E_j$ by
\begin{equation*}
  c\cdot(v_1\oplus\dots\oplus v_l)=c_1v_1\oplus\dots\oplus c_lv_l.
\end{equation*}
It follows easily from the definition that
if $\Phi$ is a slowly varying homogeneous map of $E$,
then all of the summands $\Phi_{j,\alpha}$ in the decomposition 
(\ref{flag-decomp}) are slowly varying.

A selfmap $\Phi$ of $E$ is said to be a \emph{flag map} 
(with respect to the flag $V_{\bullet}$ associated to $\oplus E_j$) 
if $\mathrm{pr}_j\Phi$ is a map of $V_j$ for all $j$, \ie
$$\Phi(v_1\oplus\dots\oplus v_l)
=\Phi_1(v_1)\oplus\Phi_2(v_1\oplus v_2)
\oplus\dots\oplus\Phi_l(v_1\oplus\dots\oplus v_l).$$
We shall say that a flag map $\Phi$ is \emph{pseudo-linear} 
if it is of the form
$$\Phi = A + H,$$
where $A$ is a linear map which preserves the splitting, and $H$ is a
polynomial flag map with no constant or linear part, such that
$\mathrm{pr}_jH=\mathrm{pr}_jH|V_{j-1}$ for $2\le j\le l$.
In other words, $\Phi$ can be written
$$\Phi(v_1\oplus\dots\oplus v_l)
=Av_1\oplus(Av_2+H_2(v_1))
\oplus\dots\oplus(Av_l+H_l(v_1\oplus\dots\oplus v_{l-1})).$$
Notice that if $A$ is invertible, then so is $\Phi$, and its inverse is
a pseudo-linear flag map whose degree is bounded in terms 
of the degree of $\Phi$.

\subsection*{Contracting pseudo-linear maps}
A bundle map $\Phi/\phi$ of a bundle $E\to X$ 
is said to be contracting if 
its linear part $d\Phi|O_E$ is contracting. If $\Phi$ is slowly varying,
then this implies that for $v\in E$ with $|v|$ small one 
has $|\Phi^Nv|\to0$ as $N\to\infty$. On the other hand, if $\Phi$ is
linear (and contracting), then this convergence holds locally uniformly for 
$v\in E$. The next result, which will be a crucial ingredient in
the proof of Theorem~\ref{pesin-u} shows that the same property 
carries over to some pseudo-linear maps $\Phi$.
\begin{prop}\label{normal-forms}
  Let $E\to X$ be a measurable bundle with a slowly varying splitting
  $\oplus E_i$ and associated flag $V_\bullet$.
  Let $\Phi/f$ be a slowly varying, pseudo-linear, contracting bundle map 
  with respect to $V_{\bullet}$. Then $\Phi^N\to0$ locally
  uniformly in the topology on $E$. 
  More precisely, there exists $\lambda<0$ with the following
  property: if $p\in X$
  and $J\subset\subset E_p$ then there exists $C=C(J,p)>0$ such that 
  \begin{equation}\label{Me5}
    |\Phi^Nv|\le Ce^{\lambda N}\quad\text{for all $v\in J$ and $N\ge0$}.
  \end{equation}
\end{prop}
\begin{proof} 
  The statement of the proposition does not change if we replace 
  the metric on $E$ by the metric given by~\eqref{adapt1},
  so let us work with that metric.
  By compactness of $J$ it suffices to show the 
  estimate~(\ref{Me5}) for large $N$.
  
  Let $\Phi = A+H$ as in the definition of pseudo-linear, and let 
  $\lambda<0$ be the associated exponent of contraction.  
  Thus we have $\Vert A\Vert\le e^\lambda$ on $E_{f^N{p}}$ for all $N\ge 0$. 
  Pick $\epsilon>0$ so small that $\lambda+3\epsilon<0$.
  
  For $N\ge1$ let $\Phi^N=\Phi^N_1\oplus\dots\oplus\Phi^N_l$ be the 
  decomposition of $\Phi^N:E\to E$ relative to the splitting 
  $E=E_1\oplus\dots\oplus E_l$. Then $\Phi^N_1=A^N$ so
  $\Vert\Phi^N_1\Vert\le e^{\lambda N}$. 
  For $2\le j\le l$ we have
  \begin{equation}\label{Me3}
    \Phi^{N+1}_j=A\Phi^N_j+H_j(\Phi^N_1\oplus\dots\oplus\Phi^N_{j-1}),
  \end{equation}
  where the $H_j$'s are polynomials with no constant or linear terms. 
  The slowly varying nature of the $H_j$'s implies the existence
  of constants $r_0>0$ and $C_0>0$ such that
  \begin{equation}\label{Me1}
    |H_j(w)|\le C_0e^{\epsilon N}|w|^2\quad
    \text{whenever $w\in E_{f^Np}$ and $|w|\le r_0e^{-\epsilon N}$}.
  \end{equation}
  We now inductively show the following estimate, which
  clearly implies the statement of the proposition
  (with $\lambda$ replaced by $\lambda+\epsilon$): there exists
  $N_0>0$ such that if $v\in J$, $1\le i\le l$ and $N\ge N_0$ then
  \begin{equation}\label{Me2}
    |\Phi^N_iv|\le e^{(\lambda+\epsilon)N}\ll r_0e^{-\epsilon N}.
  \end{equation}
  This estimate clearly holds for $i=1$.
  Suppose it holds for $1\le i<j\le l$ and let us show that it then 
  hold for $i=j$ after possibly increasing $N_0$. Indeed, if
  $N\ge N_0$ then~\eqref{Me3},~\eqref{Me1} and~\eqref{Me2} imply that
  \begin{equation}\label{Me4}
    |\Phi^{N+1}_jv|\le e^\lambda e^{(\lambda+\epsilon)N}
    +C_0e^{\epsilon N}e^{2(\lambda+\epsilon)N}
    \le e^{(\lambda+\epsilon)(N+1)}
  \end{equation}
  if $N_0$ is large enough. Thus~(\ref{Me2}) holds, which 
  completes the proof.
\end{proof}
                                
\subsection*{Regular bundle maps}
Consider a bundle $E\to X$ with a
splitting $\oplus E_j$. A linear bundle map $A:E\to E$ 
is said to be \emph{regular} if it preserves
the splitting and if there exist $\lambda_l<\dots<\lambda_1$ such that
\begin{equation*}
  \lim_{N\to\infty}\frac1N\log|A^Nv|=\lambda_j 
  \quad\text{for $v\in E_j$}
\end{equation*}
Now suppose, in addition, that the map $A$ and the splitting $\oplus E_i$
are slowly varying. 
Then one may show, using standard techniques as in~\cite{mane}, 
that there exists a slow variation $R_{\epsilon}:X\to(1,\infty)$ 
such that for all $p\in X$, all $\epsilon>0$, and all $N\in\Z$:
\begin{equation}\label{Me22}
  R_{\epsilon}(p)^{-1}e^{-\epsilon N}
  \le\frac{|A^N_pv|}{e^{\lambda_jN}|v|}
  \le R_{\epsilon}(p)e^{\epsilon N}
  \quad\text{whenever}\quad v\in E_{j,p}
\end{equation}

A general bundle map $\Phi/\phi$ is said to be regular if it is slowly
varying and if the linear bundle map $A:=dF|{0_E}$ is regular.
The vector $\underline{\lambda}=(\lambda_r,\dots,\lambda_1)$ 
is said to be the \emph{Lyapunov data} associated to $\Phi$.
It is easy to see that a regular bundle map is contracting if
and only if its Lyapunov data satisfies $\lambda_1<0$.

The part of Oseledec's Theorem that we will use 
can be stated as follows. There exist
invariant Borel sets $\klm$ in $M$,
the union of which has total measure. Over each $\klm$ there is
a measurable bundle $E^s\subset TM$ with a slowly varying 
splitting $\oplus E^{\lambda_i}$ such that $A:=df|E^s$ is a 
regular contracting bundle map fibered over $f$.

The definition of a linear bundle map being regular is asymptotic 
in nature. The concept of an \emph{adapted metric} makes this 
information easier to work with. For $\epsilon>0$ we define
\begin{equation}\label{adapt2}
  \left\langle v,w\right\rangle^*_\epsilon 
  :=\sum_{N\in\Z}
  \frac{\left\langle A^Nv,A^Nw\right\rangle}{e^{2\lambda_jN+2\epsilon|N|}}, 
  \quad v,w\in E_j
\end{equation}
and $\langle v,w\rangle_\epsilon=0$ for $v\in E_i$, $w\in E_j$, $i\ne j$.
We denote the associated norm by $|\cdot |^*_\epsilon$ and 
operator norm by $\Vert\cdot\Vert^*_\epsilon$. Often we will drop 
the subscript and write $|\cdot|^*$ instead of $|\cdot|^*_\epsilon$.

Using~\eqref{Me22} and the slow variation of the splitting $\oplus E_j$
it is not too hard to verify that
the series~\eqref{adapt2} converges, and that the metric thus
obtained is a Borel metric on $E$ with the following properties:
\begin{equation*}
  e^{\lambda_j-\epsilon}|v|^*_\epsilon
  \le |Av|^*_\epsilon 
  \le e^{\lambda_j+\epsilon}|v|^*_\epsilon
  \quad\text{for $v\in E_j$},
\end{equation*}
and there exists a slow variation 
$C_\epsilon:X\to(1,\infty)$ such that
\begin{equation}\label{Me13}
  |v|\le |v|^*_\epsilon\le C_\epsilon(p)|v|, \quad v\in E_p.
\end{equation}

Recall that a bundle map $\Phi/\phi$ is slowly varying if and only if the function
$p\mapsto\sup_m\Vert\Phi_{m,p}\Vert^{1/m}$ is slowly varying, where 
$\sum\Phi_m$ is the homogeneous expansion of $\Phi$. Using~\eqref{Me13}
it is straightforward to verify that $\Phi$ is slowly varying if and
only if there exists a slow variation $R_\epsilon$ such that
\begin{equation*}
  \Vert\Phi_{m.p}\Vert^*_\epsilon\le R_{2\epsilon}(p)^m.
\end{equation*}
In this sense, working with the given metric $|\cdot|$ or the
adapted metric $|\cdot|^*_\epsilon$ has no effect on the definition
of slow variation.

\section{Local analytic conjugation}\label{loc-conj}
The goal of this section is to establish the main step in the proof of
Theorem~\ref{pesin-u}, namely a local conjugation of regular bundle
maps to a pseudo-linear model.
\begin{thm}\label{p-and-t}
  Let $E\to X$ be a measurable bundle over an $f$-invariant
  Borel set $X\subset M$ with a slowly varying splitting $E=\oplus E_j$.
  Let $F/f$ be a slowly varying, contracting, and regular bundle map 
  of $E$ (with respect to this splitting). 
  Then there exists a polynomial bundle map $P/f$ with $dP|0_E=dF|0_E$
  and a bundle map $T/\mathrm{id}$ with $dT|0_E=\mathrm{id}$, so that
  $$T F T^{-1} = P$$
  holds on a slowly varying neighborhood of the zero section
  $O_E$ in $E$.
  Moreover, $P$ and $T$ are slowly varying, and $P$ is a pseudo-linear,
  contracting, polynomial automorphism of the flag 
  associated to the splitting of $E$.
\end{thm}

One should view $T$ as a conjugacy (\ie change of coordinates)
that conjugates $F$ to the simpler mapping $P$.
The proof of Theorem~\ref{p-and-t} 
occupies the remainder of this section.
\begin{proof}[Proof of Theorem~\ref{p-and-t}.]
  Write $A=dF|O_E$ and recall our notation $O(m)$.  
  We construct a pair of sequences $\{ P_m/f \}_{m \ge 1}$
  and $\{ T_m/\mathrm{id} \}_{m \ge 1}$ of slowly varying polynomial bundle maps
  such that:
  \begin{itemize}
  \item[(a)] $T_m = \mathrm{id} + O(2)$ and $P_m = A + O(2)$;
  \item[(b)] $P^{-1}_m T_m F-T_m = O(m+1)$;
  \item[(c)] there exists $m_0 \ge 1$ such that $P_m = P_{m_0}$ for all
    $m \ge m_0$;
  \item[(d)] $P:= P_{m_0}$ is a pseudo-linear map of the flag associated to the
    splitting of $E$;
  \item[(e)] $T_m$ converges to an analytic map $T$.
  \end{itemize}
  The construction is inductive, and proceeds as follows.  Set
  $P_1 = A$ and $T_1 = \mathrm{id}$.  Suppose that we have constructed
  $T_m$ and $P_m$. Let 
  \begin{equation*}
    T_{m+1}=T_m+H_{m+1}
    \quad\text{and}\quad
    P_{m+1}=P_m(\mathrm{id}+Q_{m+1}),
  \end{equation*}
  with $H_{m+1}$ and $Q_{m+1}$ to be determined.  
  Then a simple calculation shows that
  $$P^{-1}_{m+1}T_{m+1}F - T_{m+1} = P^{-1}_m T_m F - T_m - (Q_{m+1}+H_{m+1} -
  A^{-1}H_{m+1} A) + O(m+2).$$ 
  Thus, writing $\Phi_{m+1}:=P^{-1}_mT_mF-T_m\ \mathrm{mod}\ O(m+2)$, 
  we see that (b) holds for $m+1$ if we can
  find homogeneous solutions $H_{m+1}$ and $Q_{m+1}$  for the
  equation
  \begin{equation}\label{commutator-eq}
    \Phi_{m+1}=Q_{m+1}+H_{m+1}-A^{-1}H_{m+1}A.
  \end{equation}
  The next lemma shows that solutions of~\eqref{commutator-eq}
  exist so that (c) and (d) hold as well.
  \begin{lem}\label{comm}
    If $\Phi/\mathrm{id}$ is a homogeneous, slowly varying,
    polynomial bundle mapping of
    degree $m\ge2$, then there exist homogeneous, slowly varying
    polynomial bundle
    mappings $Q/\mathrm{id}$ and $H/\mathrm{id}$, also of degree $m$,
    such that
    \begin{equation}\label{Me10}
      \Phi= Q + H - A^{-1} H A.
    \end{equation}
    Moreover $Q$ can be chosen as follows.  
    If $m$ is sufficiently large, then one can take $Q=0$.  
    Otherwise, $Q$ is a 
    pseudo-linear map of the flag associated to $\oplus E_j$
    with no linear part, \ie $Q$ has the form
    \begin{equation*}
      Q(v_1\oplus\dots\oplus v_l)
      =0\oplus Q_2(v_1)\oplus\dots\oplus 
      Q_l(v_1\oplus\dots\oplus v_{l-1}).
    \end{equation*}
  \end{lem}
  
%  \noi \rmk The maps $Q$ obtained is this lemma will be used to make up the 
%  pseudo-linear map $P$ in Theorem~\ref{p-and-t}  (cf. Section~\ref{slow} for 
%  the definition of pseudo-linear map).  
  
  Before proving Lemma~\ref{comm}, we need to develop a few ideas.
  First, the homogeneous polynomial bundle mapping $\Phi$ 
  can be decomposed with respect to the splitting $\oplus E_i$ 
  into $(j,\alpha)$-homogeneous parts 
  $\Phi_{j,\alpha}$ as in~\eqref{flag-decomp}. 
  By linearity we only need to solve~\eqref{Me10} for each 
  summand $\Phi_{j,\alpha}$. Recall that if $\Phi$ is 
  slowly varying, then so are these summands.

  In order to solve~\eqref{Me10} we will make crucial use of the fact 
  that $A$ is regular with respect to the splitting $\oplus E_i$. 
  Let $\ul=(\lambda_1,\dots,\lambda_l)$ 
  be the Lyapunov data associated to $A$ and notice that
  \begin{equation*}
    \lambda_l<\dots<\lambda_1<0
  \end{equation*}
  since $A$ is contracting.
  If $\Phi$ is $(j,\alpha)$-homogeneous, then it follows that
  \begin{equation}\label{poly-asymp}
    \Vert A^{-N}\circ\Phi\circ A^N\Vert
    \sim e^{(\ul\cdot\alpha-\lambda_j)N}\Vert\Phi\Vert
    \quad\text{for $|N|\gg1$}.
  \end{equation}
  This estimate, which is the key to the analysis in this section,
  will be made more precise below. 
  Here we only note that the quantity $\lambda_j-\ul\cdot{\alpha}$ 
  plays an obvious role in~\eqref{poly-asymp} and leads to the 
  following definition.
  We say that the pair $(j,\alpha)$ is \emph{resonant} if 
  $$\lambda_j=\ul\cdot{\alpha}.$$
  It is called \emph{non-resonant} otherwise, 
  and more specifically,
  \emph{super-resonant} or \emph{sub-resonant} if we have 
  $\lambda_j-\ul\cdot{\alpha}<0$ or $\lambda_j-\ul\cdot{\alpha}>0$,
  respectively.
  \begin{lem}\label{b}
    There exists $m_0\in\N$ such that all pairs $(j,\alpha)$ with
    $|\alpha|>m_0$ are sub-resonant. Furthermore, if $(j,\alpha)$
    is resonant then $\alpha_j=\alpha_{j+1}=\dots=\alpha_l=0$.
  \end{lem}
  \begin{proof} 
    For the first statement pick 
    $m_0\ge\lambda_l/\lambda_1$.
    The second statement is easily verified.
  \end{proof}
  For the rest of this section we work with a fixed but small
  $\epsilon>0$. Specifically we require that
  \begin{equation}\label{Me14}
    \sup\{\ul\cdot\alpha-\lambda_j+(|\alpha|+2)\epsilon\}<0,
  \end{equation}    
  where the supremum is taken over all sub-resonant 
  pairs $(j,\alpha)$ (it is easy to check that this is possible), and
  \begin{equation}\label{Me15}
    \min\{\ul\cdot\alpha-\lambda_j-(|\alpha|+2)\epsilon\}>0
  \end{equation}    
  where the minimum is taken over all sub-resonant 
  pairs $(j,\alpha)$; the latter form a finite set by Lemma~\ref{b}. 
  
  \begin{lem}\label{a}
    With notation as in Lemma~\ref{comm}, 
    if $(j,\alpha)$ is non-resonant and $\Phi$ is 
    $(j,\alpha)$-homogeneous, then one can find a slowly
    varying $H/\mathrm{id}$ so that
    $$\Phi=H-A^{-1}HA.$$
    Further, there exists a constant $C$, depending only on
    on the Lyapunov data $\ul$, such that 
    \begin{equation}\label{Me16}
      \Vert H\Vert^*\le C\Vert\Phi\Vert^*,
    \end{equation}
    where $\Vert\cdot\Vert^*=\Vert\cdot\Vert^*_\epsilon$ is the operator
    norm associated to the adapted metric~\eqref{adapt2}.
  \end{lem}
  \begin{proof}
    The operator $\Psi\mapsto\Psi-A^{-1}\Psi A$ has the following two
    formal inverses:
    \begin{equation}\label{Me11}
      \Psi\mapsto\cg^+(\Psi):=\sum_{N=0}^{\infty}A^{-N}\Psi A^N
      \qquad\text{and}\qquad
      \Psi\mapsto\cg^-(\Psi):=-\sum_{N=1}^{\infty}A^N\Psi A^{-N}.
    \end{equation}
    We proceed in two cases.
    \begin{itemize}
    \item[(i)]  Sub-resonant case: let $H=\cg^+(\Phi)$.
    \item[(ii)] Super-resonant case: let $H=\cg^-(\Phi)$.
    \end{itemize}
    We have to show that this makes sense. Let us
    consider the sub-resonant case~(i). We have
    \begin{align*}
      \Vert\cg^+\Phi\Vert^*
      &\le\sum_{N=0}^{\infty}\Vert A^{-N}\Phi A^N\Vert^*\\
      &\le\sum_{N=0}^{\infty}
      e^{(\ul\cdot\alpha-\lambda_j+(|\alpha|+2)\epsilon )N}
      \Vert\Phi\Vert^* \\
      &=:C\Vert\Phi\Vert^*.
    \end{align*}
    It follows from~\eqref{Me14} that $C<\infty$ and that $C$ does 
    not depend on $j$, $\alpha$ or $\epsilon$. Thus~\eqref{Me16} holds,
    and this implies that $H$ is slowly varying. Notice that
    the construction of $H$ does not depend on the choice of 
    $\epsilon$.
    The super-resonant case is treated similarly, using~\eqref{Me15}.
    This completes the proof.
  \end{proof}
  \begin{proof}[Proof of Lemma~\ref{comm}.] Simply decompose $\Phi$ 
    as a sum of resonant, sub-resonant and super-resonant
    $(j,\alpha)$-homogeneous polynomial 
    bundle mappings $\Phi_{j,\alpha}$. For the non-resonant terms we
    get $H_{j,\alpha}$ from Lemma~\ref{a} and set $Q_{j,\alpha}=0$. 
    For the resonant terms we set $H_{j,\alpha}=0$ and
    $Q_{j,\alpha}=\Phi_{j,\alpha}$. Then set 
    $H=\sum H_{j,\alpha}$ and $Q=\sum Q_{j,\alpha}$. From Lemma~\ref{a}
    we get that $H$ and $Q$ are slowly varying. The remaining statements
    follow from Lemma~\ref{b}.
  \end{proof}
  
  Returning to the proof of Theorem~\ref{p-and-t} 
  we have constructed our sequences of maps $P_m$ and $T_m$, and
  have shown that they satisfy (a)--(d) above.  It remains to show that
  $T = \lim T_m$ is analytic and slowly varying. 
  For this it suffices to show that 
  $\Vert H_m\Vert^*\le B^{h(m)}$
  for some $\epsilon$-slowly varying $B:X\to(1,\infty)$ 
  and some concave function $h:\N\to\R_+$.

  Recall that the degree of the inverse of a pseudo-linear map is
  bounded in terms of the degree of the map itself. In particular,
  the degree $\nu$ of $P^{-1}$ is finite, so we may write
  $$P^{-1} = S A^{-1} = (S_1 + S_2 + \dots + S_{\nu})A^{-1},$$
  where $S_1= \mathrm{id}$ and $S_j$ is homogeneous of degree $j$, and
  all the maps are slowly varying. Finally, write $F= A G$, so that
  $G/\mathrm{id}$ is analytic and slowly varying. 
  We denote the homogeneous expansion of $G$ by
  $$G = \sum_{j=1}^{\infty} G_j,$$
  observing that $G_1 = \mathrm{id}$.
  
  Keeping $\epsilon$ fixed we now pick an $\epsilon$-slowly 
  varying $B:X\to(2,\infty)$ with the following properties: 
  \begin{align}
    \Vert S_i\Vert^*&\le B^\alpha\ \text{for $2\le i\le\nu$}\label{Me17}\\ 
    \Vert H_m\Vert^*&\le B^\alpha\ \text{for $2\le m\le m_1$}\label{Me18}\\
    \Vert G_k\Vert^*&\le B^{k\alpha}\ \text{for $k\ge 2$}\label{Me19}\\
    \Vert A^{-1}\Vert^*&\le B^\alpha
    \quad\text{and}\quad\Vert A\Vert^*\le 1.\label{Me20}
  \end{align}
  Here $\alpha\in(0,1)$ and $m_1>m_0$ will be chosen later. As before
  we write $\Vert\cdot\Vert^*=\Vert\cdot\Vert^*_\epsilon$.
  That~\eqref{Me17}--\eqref{Me20} are possible follows from the fact that
  each $S_i$ and each $H_m$ are slowly varying, as are $G$ and 
  $A^{-1}$, while $A$ is contracting.
  We claim that if $\alpha$ is small enough and $m_1$ is large enough,
  then 
  \begin{equation}\label{Me21}
    \Vert H_m\Vert^*\le B^{m-\sqrt{m}}\quad\text{for all $m\ge2$}.
  \end{equation}
  This is clearly true for $m\le m_1$ by
  our choice of $B$.  Suppose that~\eqref{Me21} holds for some $m$.
  To establish it for $m+1$ it suffices, in view of Lemma~\ref{a},
  to prove that
  \begin{equation*}
    \Vert\Phi_{m+1}\Vert^*\le C^{-1}B^{m+1-\sqrt{m+1}},
  \end{equation*}
  where $\Phi_{m+1}$ is given by~\eqref{commutator-eq}.
  To this end, let us take a look at the 
  quantity $P^{-1}T_mF-T_m$, keeping in
  mind that we are only interested in the terms of degree $m+1$.  
  Thus the terms coming from $-T_m$ automatically disappear, and we have
  \begin{align}
    \Phi_{m+1}
    &=P^{-1}T_mF-T_m\ \mathrm{mod}\ O(m+2)\notag\\
    &=\sum S_iA^{-1}H_jAG_k,\label{fmplus1}
  \end{align}
  where the sum is over all $i,j$ and $k$ 
  such that $ijk=m+1$, $k\ge1$, $1\le j\le(m+1)/2$ and $1\le i\le\nu$. 
  (The constraint on the index $j$ comes from the fact that
  (i) $1\le j\le m$, and (ii) we need only consider 
  the terms of degree $m+1$.)  

  The sum (\ref{fmplus1}) consists of considerably 
  fewer than $\nu m^2$ terms.  
  We now estimate $S_iA^{-1}H_jAG_k$ in all possible cases,
  using~\eqref{Me17}-\eqref{Me21} and the induction hypothesis,
  and assuming that $\alpha$ is small while $m_1$ is large.
  \begin{itemize}
  \item[(i)] ($k\ge2,\ i\ge 2$): then $ij\le(m+1)/2$ and
    \begin{equation*}
      \Vert S_iA^{-1}H_jAG_k\Vert^*
      \le B^{\alpha+i(\alpha+(j-\sqrt{j})+jk\alpha)}
      \le B^{(m+1)(\frac12+3\alpha)}
      \le\frac1{C\nu m^2}B^{m+1-\sqrt{m+1}}.
    \end{equation*}    
  \item[(ii)] ($k\ge2,\ i=1$): then $j\le(m+1)/2$ and
    \begin{equation*}
      \Vert A^{-1}H_jAG_k\Vert^* 
      \le B^{\alpha+(j-\sqrt{j})+jk\alpha}
      \le B^{(m+1)(\frac12+2\alpha)}
      \le\frac1{C\nu m^2}B^{m+1-\sqrt{m+1}}. 
    \end{equation*}
  \item[(iii)] ($k=1$): then $i\ge2$ and $ij=(m+1)$. Further, 
    a straightforward calculation readily shows that
    $\sqrt{m+1}(\sqrt{i}-1)-\alpha(1+i)\ge\alpha\sqrt{m+1}$, so
    \begin{align*}
      \Vert S_iA^{-1}H_jA\Vert^* 
      \le B^{\alpha+i(\alpha+(j-\sqrt{j}))}
      &=B^{(m+1-\sqrt{m+1})+(1+i)\alpha-\sqrt{m+1}(\sqrt{i}-1)}\\
      &\le B^{(m+1-\sqrt{m+1})-\alpha\sqrt{m+1}}\\
      &\le\frac1{C\nu m^2}B^{m+1-\sqrt{m+1}}.
    \end{align*}
  \end{itemize}
  Putting all this together, we have:
  \begin{equation*}
    \Vert\Phi_{m+1}\Vert^*
    \le\sum_{i,j,k}\left\Vert S_iA^{-1}H_jAG_k\right\Vert^*
    <\nu m^2\frac1{C\nu m^2}B^{m+1-\sqrt{m+1}}
    =\frac1CB^{m+1-\sqrt{m+1}}. 
    \end{equation*}
    Thus the induction step is complete and~\eqref{Me21} holds
    for all $m\ge2$. Theorem~\ref{p-and-t} now follows 
    from the definition of a slowly varying bundle map.
\end{proof}

\section{Proof of Theorem~\ref{pesin-u}}\label{2pf}
We may replace $\cK(f)$ by the Borel set $\klm$. Let us apply
Theorem~\ref{p-and-t} to the bundle $E^s\to\klm$ 
(with its Lyapunov splitting) and the bundle map 
$F=\chi^{-1}\circ f\circ\chi$ of $E^s$, where $\chi$ is the exponential
map given by Theorem~\ref{T1}. By Theorem~\ref{T2} and Theorem~\ref{T1}
the map $F$ is slowly varying, contracting and regular
(with respect to the Lyapunov splitting). Thus we 
may apply Theorem~\ref{p-and-t} to $F$.

With $P$ and $T$ as in Theorem~\ref{p-and-t} let 
$$\Psi_N:=P^{-N}T\chi^{-1}f^N.$$ 
Then for each $p\in \klm $ and each compact subset 
$J\subset\subset W_p^s$ there exists an 
integer $N_0=N_0(J)$ such that $\Psi_N$ is well
defined on $J$ whenever $N\ge N_J$.  Indeed, since $f^N$ decays 
exponentially as $N \to \infty$, and $T\chi^{-1}$ is slowly varying, $f^N$ 
will carry $J$ into the (slowly varying) domain of $T \chi^{-1}$, provided $
N$ is large enough.  Then one has
$$\Psi_{N+1} - \Psi_N = P^{-N}(P^{-1}TF-T)\chi^{-1} f^N = 0.$$
It follows that, locally uniformly in the discrete topology 
on $\cw^s$, $\Psi_N$ converges as $N\to\infty$ 
to a map $\Psi :\cw^s \to E^s$. Evidently $\Psi$
is holomorphic on the fibers $W^s_p$ of $\cw^s$, maps 
$W^s_p$ into $E^s_p$, and satisfies the functional equation
\begin{equation*}
  P^{-1}\Psi f=\Psi.
\end{equation*}
It remains only to show that $\Psi$ is bijective, which it obviously suffices
to check on fibers.  Thus we fix from here on a point $p \in \klm$.

First note that, by Proposition~\ref{inj-tube}, there exists a slowly varying
function $h$ such that $T\chi^{-1} | E^s (h)$ 
(and thus $P^{-N} T\chi^{-1}| E^s (h)$) 
is injective on fibers.  If $x_1, x_2 \in W^s_p$ are two points, 
then for large enough $N$, $f^Nx_j\in E^s(h)$ for $j=1,2$.  
(Again, this is so because $f^N$ decays exponentially, 
and thus $f^Nx_j$ eventually enters and
remains inside the slowly varying tube $E^s(h)$.)  But then, if $N$ is large 
enough,  
$$\Psi  x_1 = \Psi x_2 \Rightarrow P^{-N} T \chi^{-1} (f^Nx_1) = 
P^{-N}T\chi^{-1} (f^Nx_2) \Rightarrow f^N x_1 = f^N x_2 \Rightarrow x_1=x_2.$$
Thus $\Psi$ is injective.

Next, fix $v\in E^s_p$. We want to find $x\in W^s_p$ 
with $\Psi(x)=v$.
By Proposition~\ref{image-tube}, there exists a slowly varying tube 
$E^s(g)$ in the range of $\Psi$. 
Further, by Proposition~\ref{normal-forms}, $|P^Nv|$ decays exponentially,
so $P^N(v)$ enters  and remains inside $E^s(g)$ for large enough $N$. 
Let $y\in W^s_{f^Np}$ be such that 
$\Psi(y)=P^N(v)$ and set $x=f^{-N}y$. Then
$$\Psi(x)=\Psi(f^{-N}y)=P^{-N}\Psi(y)=v,$$
and thus $\Psi$ is surjective.

Finally, from the fact that $(dP_p)_p=A=df_p|E^s$ and that 
$(d\chi_p)_p=\mathrm{id}$, 
we see from the definition of $\Psi=\Psi_N$ that 
$(d\Psi_p)_p=\mathrm{id}$.
This completes the proof of Theorem~\ref{pesin-u}.\qed

\subsection*{Final remarks} The proof of the existence of 
the normal form in the stationary case as presented in~\cite{st,rr} 
does not use the convergence
of the formal solution $T$ of the local conjugation problem.
Instead, one shows the existence of a sufficiently large integer $m$ 
such that $P^{-N}T_mF^N$ converges. This relies heavily on the
fact that the iterates of a pseudo-linear map grow at most 
exponentially (see Lemma~1 in the Appendix of~\cite{rr}). 
This exponential growth
seems quite hard to establish in the slowly varying setting---it
would be interesting to know whether it is indeed possible. 
On the other hand, our method also provides a new proof in the
analytic stationary case.

\end{document}